\documentclass[a4paper,12pt]{article}
\usepackage{amsmath,amsthm,verbatim,amssymb,amsfonts,amscd, graphicx}
\usepackage[all]{xy}
\usepackage{xcolor}
\usepackage{geometry}
\usepackage{graphicx}
\usepackage{hyperref}
\usepackage{enumitem}
\definecolor{darkblue}{rgb}{0,0,0.3}
\definecolor{urlblue}{rgb}{0,0,0.7}
\hypersetup{
	colorlinks=true,
	citecolor=blue,
	linkcolor=darkblue,
	urlcolor=urlblue,
}

\usepackage{wrapfig}
\usepackage[pscoord]{eso-pic}% The zero point of the coordinate systemis the lower left corner of the page (the default).
\usepackage{caption}

\newcommand{\RR}{\mathbb{R}}
\newcommand{\ZZ}{\mathbb{Z}}

\newcommand{\RP}{\mathbb{RP}}

\let\div\undefined

\DeclareMathOperator{\div}{div}

\DeclareMathOperator{\sys}{sys}

\newcommand{\D}{\nabla}
\newcommand{\p}{\partial}

\renewcommand{\bar}{\overline}
\renewcommand{\tilde}{\widetilde}
\renewcommand{\epsilon}{\varepsilon}
\renewcommand{\leq}{\leqslant}
\renewcommand{\geq}{\geqslant}

\newcommand{\updatetag}[2]{}

\newtheorem{theorem}{Theorem}[section]

\newtheorem{cor}[theorem]{Corollary}

\newtheorem{defn}[theorem]{Definition}
\numberwithin{equation}{section}
\theoremstyle{definition}

\theoremstyle{definition}
\newtheorem{remark}[theorem]{Remark}
\theoremstyle{definition}

\newcommand{\tSigma}{{\tilde{\Sigma}}}
\newcommand{\tM}{{\tilde{M}}}

\newcommand{\inv}{w}

\begin{document}
	
%\title{Inverse mean curvature flow and a topological gap theorem for the $\pi_2$\,--\,systole}
\title{A topological gap theorem for the $\pi_2$\,--\,systole of positive scalar curvature 3-manifolds}
\author{Kai Xu}
\date{}
\maketitle

\begin{abstract}
	Let $M$ be a closed orientable 3-manifold with scalar curvature greater than or equal to 1. If $M$ has nonvanishing second homotopy group, then it is known that the $\pi_2$\,--\,systole of $M$ (i.e. the minimal achievable area of homotopically nontrivial spheres) is at most $8\pi$. We prove the following gap theorem: if $M$ is further not a quotient of $S^2\times S^1$, then the $\pi_2$\,--\,systole of $M$ is no greater than an improved constant $c\approx 5.44\pi$. This statement follows as a new topological application of Huisken and Ilmanen's weak inverse mean curvature flow.
\end{abstract}
	
\section{Introduction}

% MSC 2020: 53C20, 53C23, 53E10

\noindent\textbf{1.1. The main result.} Let $(M,g)$ be a closed Riemannian 3-manifold, which is assumed to be connected and oriented throughout this paper. If $M$ has nonvanishing second homotopy group, then we define the $\pi_2$\,--\,\textit{systole} of $M$ by
\[\sys\pi_2(M,g)=\inf\Big\{|S^2|_{f^*g}: f:S^2\to M\text{ is an immersion with }[f]\ne 0\in\pi_2(M)\Big\},\]
where $\pi_2(M)$ denotes the set of free homotopy classes of maps $S^2\to M$. By a theorem of Meeks and Yau \cite{Meeks-Yau_1980}, the $\pi_2$\,--\,systole is always achieved by a smooth minimizer (which is either an embedded sphere or a two-fold embedded $\RP^2$).

The interplay between curvature and geometric inequalities is one central topic in differential geometry. In \cite{Bray-Brendle-Neves_2010}, Bray, Brendle and Neves considered $\pi_2$-systolic inequalities in the context of uniformly positive scalar curvature. For manifolds $M$ with nonvanishing second homotopy group, one has the sharp inequality
\begin{equation}\label{eq-intro:BBN}
	\sys\pi_2(M,g)\cdot\min_M R_g\leq 8\pi,
\end{equation}
where $R_g$ denotes the scalar curvature of $g$. For the rigidity case, it is proved in \cite{Bray-Brendle-Neves_2010} that
\begin{equation}\label{eq-intro:rigidity}
	\sys\pi_2\cdot\min(R)=8\pi\quad\Leftrightarrow\quad\text{$M$ is isometrically covered by a round $S^2\times S^1$.}
\end{equation}
Here, a round $S^2\times S^1$ means the product of an $S^1$ with an $S^2$ with constant curvature. Inequality \eqref{eq-intro:BBN} follows from the second variation formula for the area minimizer. Note that the $\pi_2$-systole is controlled solely by scalar curvature lower bound. This is relevant with the boundedness of the Urysohn 1-width of $M$ (which suggests that $M$ is roughly a one-dimensional object), see \cite{Gromov-Lawson_1983, Liokumovich-Maximo_2023}.

Analogues of \eqref{eq-intro:BBN} were obtained in various other contexts. In Bray-Brendle-Eichmair-Neves \cite{Bray-Brendle-Eichmair-Neves_2010}, the ``$\RP^2$-systole'' of a closed 3-manifold is considered, which is defined by $\mathcal{A}:=\inf\big\{$area of embedded $\RP^2\big\}$. There it is proved that $\mathcal{A}\cdot\min(R)\leq12\pi$, with the model case being a round $\RP^3$. The non-compact version of \eqref{eq-intro:BBN} was proved by Zhu \cite{Zhu_2020_arxiv} using Gromov's $\mu$-bubble technique \cite{Gromov_2021_four_lectures}. For the case of negative scalar curvature bounds, Nunes \cite{Nunes_2013} and Lowe-Neves \cite{Lowe-Neves_2021} proved sharp area lower bounds for high-genus minimizing surfaces. The rigidity cases there are characterized respectively by product metrics $\Sigma\times\RR$ (where $\Sigma$ has genus at least 2) and hyperbolic metrics. \\

% --Optional--
% See also \cite{Zhu_2023} for another application of relevant techniques.

The main theorem of this paper is as follows.

\begin{theorem}\label{thm-intro:sys_pi2}
	Suppose $M$ is a closed 3-manifold such that $\pi_2(M)\ne0$ and $M$ is not covered by $S^2\times S^1$. Then for any metric $g$ on $M$ we have
	\begin{equation}\label{eq-intro:sys_pi2}
		\sys\pi_2(M,g)\cdot\min_M R_g\leq  24\pi\cdot\frac{2-\sqrt2}{4-\sqrt2}\qquad(\approx 5.44\pi).
	\end{equation}
\end{theorem}

Thus, the inequality \eqref{eq-intro:BBN} is improved when we topologically exclude its rigidity case \eqref{eq-intro:rigidity}. In particular, the comparison between \eqref{eq-intro:BBN} and \eqref{eq-intro:sys_pi2} shows that a universal gap is present in the $\pi_2$\,--\,systolic inequality. This is one new insight of the present work. To supply a further geometric understanding of Theorem \ref{thm-intro:sys_pi2}, we state the following result.

\begin{theorem}\label{thm-intro:local}
	Let $M$ be diffeomorphic to the 3-sphere with $k$ balls removed ($k\geq3$). Suppose $g$ is a smooth metric on $M$ such that
	\begin{enumerate}[label={(\roman*)}, nosep]
		\item all the components of $\p M$ are stable minimal surfaces;
		\item in the interior of $M$, there is no embedded stable minimal surface representing a nontrivial element in the integral homology $H_2(M,\ZZ)$.
	\end{enumerate}
	Let $A_0$ be the minimum of the area of all connected components of $\p M$. Then
	\begin{equation}\label{eq-intro:localized}
		A_0\cdot\min_M(R_g)\leq 24\pi\cdot\frac{2-\sqrt2}{4-\sqrt2}.
	\end{equation}
\end{theorem}

\begin{figure}[h]
	\vspace{6pt}
	\captionsetup{justification=centering, margin=1cm}
	\begin{center}
		\includegraphics[scale=2.2]{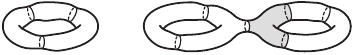}
	\end{center}
	\begin{picture}(0,0)
		%\put(78,20){$\Sigma_1$}
		\put(296,70){$\Sigma_1$}
		\put(335,86){$\Sigma_2$}
		\put(336,22){$\Sigma_3$}
	\end{picture}
	\vspace{-12pt}
	\caption{A decomposition of $S^2\times S^1$ and $(S^2\times S^1)\#(S^2\times S^1)$ by stable minimal surfaces. For the shaded piece, Theorem \ref{thm-intro:local} implies that $\min\{|\Sigma_1|,|\Sigma_2|,|\Sigma_3|\}\cdot\min(R_g)\leq(\sim 5.44\pi)$.}\label{fig:sys}
	\vspace{-6pt}
\end{figure}

Theorem \ref{thm-intro:local} suggests that Theorem \ref{thm-intro:sys_pi2} is localized at the ``non-cylindrical locus'' of the manifold, which we now demonstrate. Let $M$ have positive scalar curvature and satisfy the hypotheses of Theorem \ref{thm-intro:sys_pi2}. For simplicity, let us assume $M=(S^2\times S^1)\#(S^2\times S^1)$. By cutting along a maximal disjoint collection of stable minimal surfaces, we decompose $M$ into building pieces. This is possible for generic metrics \cite{White_2017}; similar decomposition arguments also appear in \cite{Liokumovich-Maximo_2023, Sormani_2023}. Each building piece is diffeomorphic to a sphere with disks removed, and its boundary consists of stable minimal spheres. See the right hand side of Figure \ref{fig:sys} for a depiction of the decomposition.
The topology of $M$ implies the following: there must exist a non-cylindrical piece, namely, a piece with three or more boundary components (such as the shaded piece in Figure \ref{fig:sys}). Theorem \ref{thm-intro:local} then yields the improved area bound \eqref{eq-intro:localized} for the boundary of such a piece. With a modification of the decomposition argument, we can arrange that all the boundary surfaces are nontrivial in $\pi_2(M)$. In this way, we may recover the global inequality \eqref{eq-intro:sys_pi2} from the local version \eqref{eq-intro:localized}. For the full argument that follows the idea presented here, see Remark \ref{rmk-sys:decomp}. We remind that the formal proof of the main theorems goes in the opposite direction: we first prove Theorem \ref{thm-intro:sys_pi2}, then obtain Theorem \ref{thm-intro:local} via a doubling argument.

The present work is partially motivated by the stability problem for scalar curvature. Concerning the inequality \eqref{eq-intro:BBN}, the stability problem asks the following: if we relax the rigidity condition in \eqref{eq-intro:rigidity} to almost rigidity, i.e. if we impose $\sys\pi_2\cdot\min(R)\geq 8\pi-\epsilon$, whether the manifold $M$ remains close to a round cylinder (or its quotient) in some sense. In this respect, Theorem \ref{thm-intro:sys_pi2} confirms that $M$ must be topologically the same as rigidity, thus establishing a \textit{topological stability}. On the other hand, the metric stability of \eqref{eq-intro:BBN} remains an open question. The weak control of scalar curvature on the metric allows various pathological phenomena. We refer the reader to the survey of Sormani \cite{Sormani_2023}, and references therein, for an overview of this topic.

\vspace{9pt}

\noindent\textbf{1.2. Proof by weak inverse mean curvature flow.} The method for proving Theorem \ref{thm-intro:sys_pi2} is the weak inverse mean curvature flow (abbreviated by IMCF), which was introduced by Huisken and Ilmanen in \cite{Huisken-Ilmanen_2001}. This is a suitable weak version of the parabolic flow that evolves a hypersurface at the speed of $1/H$, where $H$ denotes the mean curvature. The weak IMCF has seen remarkable success in studying positive scalar curvature in dimension three. One of its most renowned applications is the Riemannian Penrose inequality for connected horizons \cite{Huisken-Ilmanen_2001}. Further applications can be found in for example \cite{Bray-Neves_2004, Huisken-Koerber_2023, Shi_2016}. A brief introduction of this flow is contained in Section \ref{sec:prelim}.

The proof of Theorem \ref{thm-intro:sys_pi2} starts with the following setup. Let $\Sigma\subset M$ an area minimizer in $\pi_2(M)$. Consider the universal cover of $M$ (denoted by $\tM$), and a lift of $\Sigma$ onto $\tM$ (denoted by $\tSigma$). Let $\{\Sigma_t\}$ solve the weak IMCF in $\tM$, with initial value $\tSigma$. Figure \ref{fig:imcf} below illustrates the behavior of the flow in the universal cover of $(S^2\times S^1)\#(S^2\times S^1)$. In this case, $\tM$ resembles the shape of a binary tree. The initial surface $\tSigma$ sits at the bottom of the figure, and the flow moves upwards.
\begin{figure}[h]
	\captionsetup{justification=centering, margin=1.5cm}
	\begin{center}
		\includegraphics[scale=2]{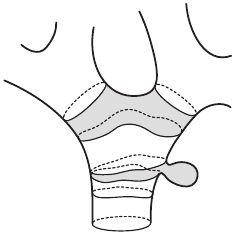}
	\end{center}
	\vspace{-12pt}
	\caption{The weak inverse mean curvature flow on $\tM$ starting with $\tSigma$. The grey regions represent ``jumpings'' in the weak flow.}\label{fig:imcf}
	\vspace{12pt}
	\begin{picture}(0,0)
		\put(270,260){$\tM$}
		\put(188,67){$\tSigma$}
		\put(187,98){$\Sigma_t$}
		\put(176,141){$\Sigma_T$}
		\put(158,163){$\Sigma_T^+$}
		\put(154,205){$\geq\!\sys\pi_2$}
		\put(285,204){$\geq\!\sys\pi_2$}
		\put(262,68){$\text{Area}=\sys\pi_2$}
	\end{picture}
	\vspace{-24pt}
\end{figure}
The weak IMCF exhibits ``jumping'' behaviors, which are represented by the shadowed regions. At a certain time, as labeled by $T$ in the figure, the flow splits into two components which will individually continue evolving in their own branches.

\vspace{12pt}

We remark on several key aspects of the proof.

1. A challenge for approaching Theorem \ref{thm-intro:sys_pi2} is the need to distinguish $S^2\times S^1$ (and its quotients) from other manifolds. In the IMCF approach, the existence theory plays the role of telling apart topological types. The weak flow needs to satisfy a condition called \textit{properness} to have desirable analytic properties. The properness of the weak flow holds under the assumptions of Theorem \ref{thm-intro:sys_pi2}; see the references \cite{Mari-Rigoli-Setti_2022, Xu_2023}, Theorem \ref{thm-imcf:exist}, and the argument at the beginning of Section \ref{sec:sys}. In contrast, the universal cover of $S^2\times S^1$ does not admit any proper weak solution, and our proof no longer works in this case.

2. An important property of the weak IMCF is the \textit{Geroch monotonicity formula}
\begin{equation}\label{eq-intro:monotonicity}
	\frac d{dt}\int_{\Sigma_t}H^2\leq -\frac12\int_{\Sigma_t}H^2+4\pi\chi(\Sigma_t)-\int_{\Sigma_t}R_g \quad\text{(in a suitable weak sense)},
\end{equation}
where $H$ is the mean curvature of $\Sigma_t$. This formula lies at the heart of our proof, as well as many other applications \cite{Bray-Neves_2004, Huisken-Ilmanen_2001, Huisken-Koerber_2023, Shi_2016}. To make use of \eqref{eq-intro:monotonicity}, one needs to control the Euler characteristic term on the right hand side.

3. The control of Euler characteristic relies on certain topological properties of the weak flow. Here the relevant property is:
\begin{equation}\label{eq-intro:topology}
	\text{every spherical component of $\Sigma_t$ is nontrivial in $\pi_2(\tM)$.}
\end{equation}
This immediately implies that $\chi(\Sigma_t)\leq2$ as long as $|\Sigma_t|<2\sys\pi_2(\tM)$. In combination with the exponential growth law $|\Sigma_t|=e^t|\Sigma_0|=e^t\sys\pi_2(M)$ and the fact $\sys\pi_2(\tM)=\sys\pi_2(M)$, we obtain $\chi(\Sigma_t)\leq2$ for all $t<\log 2$, which is the desired bound on Euler characteristics. This intuitively means that the first splitting time (i.e. the value $T$ in Figure \ref{fig:imcf}) is at least $\log 2$.

Finally, the main inequality \eqref{eq-intro:sys_pi2} is obtained by ``solving'' \eqref{eq-intro:monotonicity} in the time interval $[0,\log 2]$ and applying the observation in item 3. Based on the present proof, the sharpness of \eqref{eq-intro:sys_pi2} or the case of rigidity is unknown.

\vspace{12pt}

\noindent\textbf{1.3 Further questions.} Several further questions appear in the direction of Theorem \ref{thm-intro:sys_pi2}. The topological invariant
\[\inv(M)=\sup\Big\{\sys\pi_2(M,g)\cdot\min_M R_g: g\text{ is a Riemannian metric on }M\Big\}\]
is now known to be not trivial. Indeed, we know that $w(M)=+\infty$ if $M$ is a sphere quotient, and $\inv(S^2\times S^1)=\inv(\RP^3\#\RP^3)=8\pi$, and $\inv(M)$ is no greater than roughly $5.44\pi$ for all the other $M$. It is natural to ask whether further gaps are present for other topological types. Reasonable candidates for this include $M=L(p,q)\# L(p,q)$, where $L(p,q)$ denotes a lens space. The sharpness of the constant in Theorem \ref{thm-intro:sys_pi2}, and the existence of extremal metrics, are not known as well.

\vspace{12pt}

\textbf{Acknowledgements}. The author appreciates Hubert Bray, Yupei Huang, Demetre Kazaras and Marcus Khuri for many insightful discussions and helpful comments on the previous versions of this paper, and the anonymous referee for the valuable suggestions on improvements.

% New section. ------------------------------------------------------------

\section{Preliminaries}\label{sec:prelim}

\vspace{9pt}

\ \ \ \ The main technique that we use is Huisken and Ilmanen's weak formulation of the inverse mean curvature flow (IMCF) \cite{Huisken-Ilmanen_2001}. Instead of diving into the analytic nature of weak IMCF, we reduce the technicality by only including statements that are necessary for proving Theorem \ref{thm-intro:sys_pi2}. We refer the reader to \cite{Huisken-Ilmanen_2001, Lee_rel} for comprehensive introductions of this rich subject.

The classical IMCF is a parabolic flow that starts with a closed strictly mean convex hypersurface $\Sigma_0$ and evolves the hypersurface by the inverse of its mean curvature:
\begin{equation}\label{eq-imcf:smooth_flow}
	\frac{\p\Sigma_t}{\p t}=\frac1H\nu.
\end{equation}
The IMCF can be alternatively described using a parametrizing function. Let $u$ be a function defined such that $\Sigma_t=\{u=t\}$. Then one can show that $u$ solves the degenerate-elliptic equation
\begin{equation}\label{eq-imcf:deg_elliptic}
	\div\Big(\frac{\D u}{|\D u|}\Big)=|\D u|.
\end{equation}
For a star-shaped strictly mean convex hypersurface in $\RR^n$, the classical IMCF exists for all time \cite{Gerhardt_1991, Urbas_1990}. In many cases that naturally arise, for example when the exterior region is not diffeomorphic to $\Sigma_0\times[0,\infty)$, the smooth flow forms singularity and terminates in finite time. This makes it necessary to consider a weak version of the flow. The weak formulation in \cite{Huisken-Ilmanen_2001} heuristically avoids singularity formation through a jumping mechanism, in which there are countably many times $t$ when $\Sigma_t$ immediately ``jumps'' to its outermost area-minimizing envelope.

Analytically, the weak formulation is defined by finding suitable variational principles for \eqref{eq-imcf:deg_elliptic}. The initial value problem of weak IMCF is defined as follows.

\begin{defn}[{\cite[p.367]{Huisken-Ilmanen_2001}}]\label{def-imcf:weak_sol}
	Let $M$ be a connected, complete, non-compact manifold and $E_0\subset M$ be a bounded smooth domain. A locally Lipschitz function $u$ on $M$ is called a weak solution of the IMCF with initial condition $E_0$, if
	
	(1) $E_0=\{u<0\}$,
	
	(2) For any compact set $K\subset\subset M\setminus\bar{E_0}$ and any locally Lipschitz function $v$ such that $\{u\ne v\}\subset K$, we have
	\[\int_K|\D u|+u|\D u|\leq\int_K|\D v|+v|\D u|.\]
	
	The weak solution $u$ is called proper if the set $E_t:=\{u<t\}$ is bounded for all $t\geq0$, equivalently if $\lim_{x\to\infty}u(x)=+\infty$ uniformly.
\end{defn}

Note that the specific value of $u$ in $E_0$ does not affect $u$ being a weak solution. We next record the following existence criterion and properties of weak solutions.

\begin{theorem}\label{thm-imcf:exist}
	Suppose $M$ is a closed 3-manifold whose fundamental group has exponential growth. Let $\tM$ be the universal cover of $M$. Then starting with any smooth bounded domain $E_0$ in $\tM$, there exists a unique proper weak solution of IMCF in the sense of Definition \ref{def-imcf:weak_sol}.
\end{theorem}
\begin{proof}
	By a theorem of Coulhon and Saloff-Coste \cite{Coulhon-Saloff-Coste_1993}, $\tM$ supports a uniform Euclidean isoperimetric inequality $|\p^*E|\geq c|E|^{2/3}$, $\forall E\subset\subset\tM$. Then we apply \cite[Theorem 1.2]{Xu_2023} to prove the theorem.
\end{proof}

\begin{theorem}\label{thm-imcf:properties}
	Let $M$ be a connected, complete, non-compact 3-manifold, and $u$ is a proper weak solution of IMCF on $M$ with initial condition $E_0$ (which is a bounded smooth domain). Assume that $E_0$ is connected and is outward minimizing, in the sense that any bounded domain $E\supset E_0$ satisfies $|\p^*E|\geq|\p E_0|$. Then $u$ satisfies the following properties:
	
	(1) For all $t\geq0$, $\Sigma_t:=\p E_t$ is a $C^{1,\alpha}$ hypersurface. The weak mean curvature of $\Sigma_t$ (denoted by $H$) exists for all $t$ and is equal to $|\D u|$ almost everywhere for almost every $t$.
	
	(2) The area of $\Sigma_t$ grows exponentially: $|\Sigma_t|=e^t|\Sigma_0|$.
	
	(3) For all $t>0$, $E_t$ is connected and $M\setminus E_t$ has no compact connected component.
	
	(4) Let $R$ be the scalar curvature of $M$. Then% function $H(t)=\int_{\Sigma_t}H^2$ satisfies
	\begin{equation}\label{eq-imcf:monotonicity2}
		\begin{aligned}
			\int_{\Sigma_{t_2}}H^2 &\leq \int_{\Sigma_{t_1}}H^2-\int_{t_1}^{t_2}\int_{\Sigma_s}\Big[2\frac{|\D_{\Sigma_s}H|^2}{H^2}+\frac12|\mathring{A}|^2\Big]\,ds \\
			&\qquad\qquad\qquad +\int_{t_1}^{t_2}\Big[4\pi\chi(\Sigma_s)-\int_{\Sigma_s}R-\frac12\int_{\Sigma_s}H^2\Big]\,ds,\quad \forall\,0\leq t_1<t_2,
		\end{aligned}
	\end{equation}
	where $\mathring{A}$ denotes the traceless weak second fundamental form of $\Sigma_s$.
\end{theorem}

The four properties in the theorem follow from Theorem 1.3, equation (1.12), Lemma 1.6, the proof of Lemma 4.2, and Theorem 5.7 in \cite{Huisken-Ilmanen_2001}. Properties (2) and (4) are not hard to verify for smooth flow \eqref{eq-imcf:smooth_flow}. Property (4) is usually referred to as the \textit{monotonicity of Hawking mass}, which is a key step in Huisken and Ilmanen's proof of the Riemannian Penrose inequality. Note that $H(t)$ is not necessarily continuous along the weak flow, thus property (4) must be expressed in integral form.

The following topological statement is a direct consequence of (3).

\begin{cor}\label{cor-imcf:topology}
	For every $t>0$, the map $H_2(\Sigma_t,\ZZ)\to H_2(M\setminus E_0,\ZZ)$ induced by embedding is injective.
\end{cor}
\begin{proof}
	Denote $M_t=\bar{E_t}\setminus E_0$ and $M'_t=M\setminus E_t$. By property (3) and the Lefschetz duality \cite[Theorem 3.43 and Exercise 3.3.35]{Hatcher}, we have $H_3(M_t,\Sigma_t,\ZZ)=H^0(M_t,\Sigma_0,\ZZ)=0$ and $H_3(M'_t,\Sigma_t,\ZZ)=H^0_c(M'_t,\ZZ)=0$. Hence by the excision theorem \cite[Theorem 2.20]{Hatcher},
	% Detail: see Hatcher Ex. 3.3.35 for noncompact form of Lefschetz duality.
	\[H_3(M\setminus E_0,\Sigma_t,\ZZ)\cong H_3(M_t,\Sigma_t,\ZZ)\oplus H_3(M'_t,\Sigma_t,\ZZ)=0.\]
	The corollary then follows from the long exact sequence of relative homology
	% Detail: note that $H_3(M\setminus E_0,\Sigma_t,\ZZ)$ is isomorphic to $H_3(M\setminus E_0,N_\epsilon(\Sigma_t),\ZZ)$, where $N_\epsilon$ denotes collar neighborhood. Use excision theorem for the latter.
	\[0=H_3(M\setminus E_0,\Sigma_t,\ZZ)\to H_2(\Sigma_t,\ZZ)\to H_2(M\setminus E_0,\ZZ)\to\cdots.\qedhere\]
\end{proof}

From \eqref{eq-imcf:monotonicity2} we obtain the following:

\begin{cor}\label{cor-imcf:Gronwall_ineq}
	Let $t>0$. If $R\geq\lambda$ on $E_t\setminus E_0$, then
	\begin{equation}\label{eq-imcf:monotonicity}
		\int_{\Sigma_t}H^2\leq e^{-t/2}\int_{\Sigma_0}H^2+4\pi e^{-t/2}\int_0^t e^{s/2}\chi(\Sigma_s)\,ds-\frac23\lambda|\Sigma_0|\big(e^t-e^{-t/2}\big).
	\end{equation}
\end{cor}
\begin{proof}
	For $0\leq s\leq t$, we denote $H(s)=\int_{\Sigma_s}H^2$, $\chi(s)=\chi(\Sigma_s)$, $F(s)=\int_s^t H(r)\,dr$. Applying \eqref{eq-imcf:monotonicity2} with $t_1=s$, $t_2=t$, and noting that $\int_{\Sigma_r}R\geq\lambda e^r|\Sigma_0|$, for almost every $s$ it holds
	\[\begin{aligned}
		\frac d{ds}\big[e^{s/2}F(s)\big] &= e^{s/2}\Big[\frac12F(s)-H(s)\Big] \\
		&\leq 4\pi e^{s/2}\int_s^t\chi(r)\,dr-\lambda|\Sigma_0|e^{s/2}\big(e^t-e^s\big)-e^{s/2}H(t).
	\end{aligned}\]
	Integrating from 0 to $t$, this implies
	\begin{equation}\label{eq-imcf:aux1}
		-F(0)\leq4\pi\int_0^t 2\big(e^{r/2}-1\big)\chi(r)\,dr-\lambda|\Sigma_0|\Big(\frac43e^{3t/2}-2e^t+\frac23\Big)-2\big(e^{t/2}-1\big)H(t).
	\end{equation}
	Applying \eqref{eq-imcf:monotonicity2} once more with $t_1=0$, $t_2=t$, and plugging in \eqref{eq-imcf:aux1}, we obtain \eqref{eq-imcf:monotonicity}.
\end{proof}

\section{Proof of the main theorems}\label{sec:sys}

\ \ \ \ We first prove Theorem \ref{thm-intro:sys_pi2}. Since the statement is vacuous when the scalar curvature of $M$ is non-positive somewhere, we may assume that $M$ has positive scalar curvature. It is known that closed 3-manifolds with positive scalar curvature are diffeomorphic to connected sums of spherical space forms and $S^2\times S^1$, see \cite[Theorem 1.29]{Lee_rel}. Thus we write
\begin{equation}\label{eq-sys:topology}
	M=(S^2\times S^1)\#\cdots\#(S^2\times S^1)\#(S^3/\Gamma_1)\#\cdots\#(S^3/\Gamma_k).
\end{equation}
Based on this classification, there are three possible cases of the growth of $\pi_1(M)$:
\begin{enumerate}[label=(\theenumi),topsep=2pt,itemsep=-1ex,partopsep=1ex,parsep=1ex]
	\item $M$ is a spherical space form, for which $\pi_1(M)$ is finite.
	\item $M$ is either $S^2\times S^1$ or $\RP^3\#\RP^3$, for which $\pi_1(M)$ is virtually cyclic.
	\item All the remaining cases, where $\pi_1(M)$ has exponential growth.
\end{enumerate}
Therefore, the topological condition in Theorem \ref{thm-intro:sys_pi2} implies case (3).

\begin{proof}[Proof of Theorem \ref{thm-intro:sys_pi2}] {\ }
	
	Assume that $\min_M R=\lambda>0$, and denote $A_0=\sys\pi_2(M)$. Taking a double cover if necessary, we assume that $M$ has no $\RP^3$ factors. (From the definition, the $\pi_2$\,--\,systole of a manifold is equal to that of any covering space.) Then by a theorem of Meeks and Yau \cite{Meeks-Yau_1980}, there exists an embedded sphere $\Sigma\subset M$ with $|\Sigma|=A_0$. Let $\tM$ be the universal cover of $M$, and let $\tSigma$ be an isometric lift of $\Sigma$ onto $\tM$. Note that $\tSigma$ must be separating. Otherwise, there is a loop in $\tM$ that intersects $\tSigma$ once, violating $\pi_1(\tM)=0$. Denote the two connected components of $\tM\setminus\tSigma$ by $\tM_1$ and $\tM_2$.
	
	\vspace{9pt}
	\noindent\textbf{Claim 1.} Both $\tM_1$, $\tM_2$ are non-compact.
	\begin{proof}
		It follows from van Kampen's theorem that $\pi_1(\tM_1)=0$. If $\tM_1$ is compact, then by the relative Poincare duality and the long exact sequence of relative cohomology, we have $H_2(\tM_1,\ZZ)=H^1(\tM_1,\p\tM_1,\ZZ)=H^1(\tM_1,\ZZ)=0$. Hence $\pi_2(\tM_1)=0$ by Hurewicz's theorem. This contradicts with the fact that $[\tSigma]\ne0\in\pi_2(\tM)$. For the same reason, $\tM_2$ is noncompact.
	\end{proof}
	
	\begin{figure}[h]
		\captionsetup{justification=centering, margin=1.5cm}
		\begin{center}
			\includegraphics[scale=2.5]{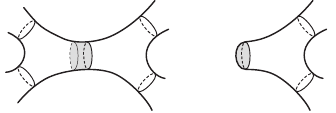}
		\end{center}
		\vspace{-6pt}
		\caption{Before and after the excision.}\label{fig:excision}
		%\vspace{12pt}
		\begin{picture}(0,0)
			\put(28,103){$\tM_1$}
			\put(215,100){$\tM_2$}
			\put(127,126){$\tSigma$}
			\put(116,76){$E'_0$}
			\put(143,103){$\to$}
			\put(152,111){\small IMCF\normalsize}
			\put(90,103){$\leftarrow$}
			\put(65,111){\small IMCF\normalsize}
			\put(275,83){$D=E_0$}
			\put(404,100){$\tM_2$}
			\put(317,126){$\tSigma$}
			\put(334,103){$\to$}
			\put(343,111){\small IMCF\normalsize}
		\end{picture}
		\vspace{-18pt}
	\end{figure}
	
	Let $E'_0\subset\tM_1$ be a (one-sided) collar neighborhood of $\tSigma$ in $\tM_1$. By Theorem \ref{thm-imcf:exist} and the fact that $\pi_1(M)$ has exponential growth, there exists a proper weak solution $u'$ of IMCF on $\tM$, with initial condition $E'_0$. As the part of the weak solution in $\tM_1$ is not utilized and causes inconvenience, we excise it from the manifold in the following way. Consider a new Riemannian manifold $N=\tM_2\cup_{\tSigma} D$ ($D$ denotes a 3-disk), where we smoothly extend the metric on $\tM_2$ into $D$. See Figure \ref{fig:excision} for how this operation is done. Let $u$ be a locally Lipschitz function on $N$ that coincides with $u'$ on $\tM_2$ and is negative in $D$. By checking Definition \ref{def-imcf:weak_sol}, one verifies that $u$ is a proper weak solution on $N$ with initial condition $E_0=D$.

	\vspace{9pt}
	\noindent\textbf{Claim 2.} $E_0$ is outward minimizing in $N$.
	\begin{proof}
		By \cite[Corollary 3.15]{Xu_2023} or \cite[Proposition 3.1]{Fogagnolo-Mazzieri_2022}, $E'_0$ admits a bounded least area envelope $F'_0$ in $\tM$, where for ``least area evenlope'' one means that $F'_0$ minimizes the perimeter among all sets $F'$ with $E'_0\subset F'\subset\subset\tM$. It follows directly from this definition that $F_0:=(F'_0\cap\tM_2)\cup E_0$ is a least area envelope of $E_0$ in $N$. The claim follows if we can show that $|\p F_0|\geq|\p E_0|$.
		%\[|\p^*F_0|=\inf\big\{|\p^*F|: E_0\subset F\subset\subset N\big\}.\]
		By the strong maximum principle, $F_0$ either coincides with $E_0$ or has a stable minimal boundary in $N\setminus\bar{E_0}\cong\tM_2$. The former case immediately implies the claim. Suppose that the latter holds. Then we have $[\p F_0]=[\tSigma]$ viewed as elements in $H_2(\tM,\ZZ)$, which is nonzero by Hurewicz theorem. Since $\tM$ has positive scalar curvature, all connected components of $\p F_0$ are spherical. Hence at least one component of $\p F_0$ is nonzero in $\pi_2(\tM)$. Finally we have $|\p F_0|\geq\sys\pi_2(\tM)=\sys\pi_2(M)=A_0=|\p E_0|$. This shows that $E_0$ is outward minimizing.
	\end{proof}
	
	Claim 2 and Theorem \ref{thm-imcf:properties} (2) then imply that $|\Sigma_t|=e^tA_0$ for all $t>0$. Define
	\[T=\inf\big\{t>0: \Sigma_t\text{ has at least two spherical connected components}\big\}.\]
	
	\noindent\textbf{Claim 3.} Each spherical component of $\Sigma_t$ ($t>0$) has area at least $A_0$.
	\begin{proof}
		Let $\Sigma'$ be a spherical component of $\Sigma_t$, which we may also view as a surface in $\tM_2$ or in $\tM$. By Corollary \ref{cor-imcf:topology}, $[\Sigma']$ is nonzero in $H_2(N\setminus E_0,\ZZ)\cong H_2(\tM_2,\ZZ)$. By the fact that $\tM_1$ is noncompact (Claim 1) and the long exact sequence of relative homology
		\[0=H_3(\tM,\tM_2,\ZZ)\to H_2(\tM_2,\ZZ)\to H_2(\tM,\ZZ),\]
		$[\Sigma']$ is nonzero in $H_2(\tM,\ZZ)$. Since $\Sigma'$ is spherical, it is a nonzero element in $\pi_2(\tM)$. Hence $|\Sigma'|\geq\sys\pi_2(\tM)=\sys\pi_2(M)=A_0$.
	\end{proof}
	
	By the definition of $T$, there is a sequence $t_i\to0$ such that $\Sigma_{T+t_i}$ has more than one spherical component. Therefore $e^{T+t_i}A_0=|\Sigma_{T+t_i}|\geq 2A_0$. Letting $i\to\infty$ it follows that $T\geq\log2$. Therefore, $\chi(\Sigma_t)\leq 2$ for all $0\leq t<\log 2$. Finally, we utilize the monotonicity formula \eqref{eq-imcf:monotonicity} to obtain
	\[0 \leq \int_{\Sigma_t}H^2
	\leq 16\pi\big(1-e^{-t/2}\big)-\frac23\lambda A_0\big(e^t-e^{-t/2}\big),\quad\forall\,t\leq\log2.\]
	Taking $t=\log 2$ we have
	\[0 \leq 16\pi\big(1-\frac1{\sqrt 2}\big)-\frac23\lambda A_0\big(2-\frac1{\sqrt 2}\big),\]
	which implies \eqref{eq-intro:sys_pi2}. This completes the proof of Theorem \ref{thm-intro:sys_pi2}.
\end{proof}

Finally, we prove Theorem \ref{thm-intro:local} using a doubling argument.

\begin{proof}[Proof of Theorem \ref{thm-intro:local}]
	We can assume $\min_M R=\lambda>0$. Let $\Sigma_i$ ($1\leq i\leq k$) be all the connected components of $\p M$. Denote $h_i=g|_{\Sigma_i}$ the restricted metrics on $\Sigma_i$, and choose $\varphi_i>0$ as the first eigenfunctions of the stability operator for $\Sigma_i$. Therefore $\varphi_i$ satisfy
	\[\Delta_{h_i}\varphi_i\leq\Big(K_{h_i}-\frac12\lambda\Big)\varphi_i,\]
	where $K$ denotes the Gauss curvature. For $T$ a sufficiently large constant to be chosen, let $P_i$ be diffeomorphic to the cylinders $\Sigma_i\times[-T,T]$ and equipped with the warped product metrics $g_i=h_i+\varphi_i^2dt^2$ ($-T\leq t\leq T$). The scalar curvature of $P_i$ is given by
	\[R_{g_i}=2\Big(K_{h_i}-\frac{\Delta_{h_i}\varphi_i}{\varphi_i}\Big)\geq\lambda,\]
	see \cite[Proposition 1.13]{Lee_rel}. Let $M^\pm$ be identical copies of $M$, whose metrics are still denoted by $g$. Set
	\[N=M^-\cup_{\p M^-}\Big(\bigsqcup_{i=1}^k P_i\Big)\cup_{\p M^+}M^+,\]
	equiped with the Lipschitz metric $g_N$ that agrees with $g$ on $M^\pm$ and agrees with $g_i$ on $P_i$. Topologically, $N$ is a connected sum of $k-1$ copies of $S^2\times S^1$.
	
	Given any $0<\epsilon<1/100$, we claim the following:
	
	\vspace{9pt}
	\noindent\textbf{Claim 1}. There is a sufficiently large $T$ for which the following holds. Suppose $g'_N$ is any smooth metric on $N$, such that $g'_N=g_i$ in each $P_i$, and $||g'_N-g_N||_{C^0(g_N)}\leq\epsilon^3$. Then $\sys\pi_2(N,g'_N)\geq(1-\epsilon)A_0$.
	
	\vspace{9pt}
	\noindent\textbf{Claim 2}. There exists a smooth metric $g'_N$ that satisfies the hypotheses of Claim 1, and moreover has $R_{g'_N}\geq\lambda-\epsilon$.
	
	\vspace{9pt}
	The smoothness of $g'_N$ in the claims makes sense as a set of coordinate charts across $\p M^\pm$ will be specified in the proof of Claim 2. With the two claims, it follows from Theorem \ref{thm-intro:sys_pi2} that $(\lambda-\epsilon)\cdot(1-\epsilon)A_0\leq c$, which by $\epsilon\to0$ proves Theorem \ref{thm-intro:local}. \qedhere
	
	\begin{proof}[Proof of Claim 1] {\ }
		
		We start with proving two technical facts.
		%if $T$ is sufficiently large in the definition of $P_i$, then an area minimizer in $\pi_2(N,g'_N)$ is either a horizontal slice in one of the $P_i$, or does not intersect any $P_i$. To prove this, we notice the following:
		
		\vspace{3pt}
		\noindent\textbf{Claim 3.} There exists a constant $c_0$ independent of $T$ such that: if $t\in(-T+2,T-2)$, $1\leq i\leq k$, and $S\subset N$ is a closed $g'_N$-minimal surface that intersects with $\Sigma_i\times(t-1,t+1)\subset P_i$, then $\big|S\cap\big(\Sigma_i\times(t-2,t+2)\big)\big|_{g_i}\geq c_0$.
		
		\vspace{3pt}
		\noindent\textit{Proof.} Let $r_0>0$ be such that $B_{g_i}(x,r_0)\subset\Sigma_i\times(t-1.5,t+1.5)$ for all $x\in\Sigma_i\times(t-1,t+1)$. By the classical monotonicity formula \cite[Lemma 1]{Meeks-Yau_1980}, there exists a constant $c>0$ such that $|S\cap B(x,r)|_{g_i}\geq cr^n$ for all $x\in S\cap\big(\Sigma_i\times(t-1,t+1)\big)$ and $r\in[0,r_0]$. The constants $c,r_0$ depend only on the geometry of $\Sigma_i\times(t-2,t+2)$, thus not on $T$ by the translation symmetry. Now $c_0=c_1r_0^n$ has the desired property.
		
		%The justification of this fact relies on the classical monotonicity formula for minimal surfaces \cite[Lemma 1]{Meeks-Yau_1980}. In fact, there exists a constant $c_0$ independent of $T$ such that: if a closed $g'_N$-minimal surface $S\subset N$ intersects $\Sigma_i\times(t-1,t+1)\subset P_i$ for some $t\in(-T+2,T-2)$ and some $i$, then $|S\cap\big(\Sigma_i\times(t-2,t+2)\big)|_{g_i}\geq c_0$. To see this, we first choose a small $r_0$ such that $B_{g_i}(x,r_0)\subset\Sigma_i\times(t-2,t+2)$ for any $x\in\Sigma_i\times(t-1,t+1)$. By the monotonicity formula, there exists $r_1\leq r_0$ and $c_1>0$ such that $|S\cap B(x,r_1)|_{g_i}\geq c_1r_1^n$ whenever $x\in S\cap\big(\Sigma_i\times(t-1,t+1)\big)$. The constants $c_1,r_0,r_1$ depend only on $\{h_i\}$ and $\{\varphi_i\}$ but not on $T$. Now $c_0=c_1r_1^n$ has the desired property.
		
		\vspace{3pt}
		\noindent\textbf{Claim 4.} Set $T=16A_0/c_0$. Then an $g'_N$-area minimizer in $\pi_2(N)$ either coincides with $\Sigma_i\times\{t\}\subset P_i$ for some $i$, or does not intersect any $P_i$.
		
		\vspace{3pt}
		\noindent\textit{Proof.} Let $S$ be such a minimizer. Suppose $S$ intersects with the interior of some $P_i$. By the strong maximum principle, $S$ either intersects with $\Sigma_i\times\{t\}$ for every $-T<t<T$, or coincides with $\Sigma_i\times\{t\}$ for some $t$. If the former case happens, then applying Claim 1 to each $\Sigma_i\times[4j,4j+4]$, $|j|\leq\lfloor T/4\rfloor$, we obtain
		\[|S|_{g'_N}\geq|S\cap\big(\Sigma_i\times(-T,T)\big)|_{g_i}\geq c_0\cdot2\lfloor\frac T4\rfloor>2A_0.\]
		This contradicts the minimality of $S$, since one of the $\Sigma_i$ has area $A_0$.
		
		\vspace{3pt}
		Back to the proof of claim 1, we choose $T$ as in Claim 4. Let $S$ be a $g'_N$-area minimizer in $\pi_2(N)$, which is a smoothly embedded sphere. Thus Claim 4 applies to $S$. If $S$ is a horizontal slice in some $P_i$, then $|S|_{g'_N}=|\Sigma_i|_{h_i}\geq A_0$, implying Claim 1. Now we may assume that $S$ does not intersect all the $P_i$.
		
		%Suppose that $S$ has nonempty intersection with the interior of some $P_i$. By the strong maximum principle, $S$ either intersects with $\Sigma_i\times\{t\}\subset P_i$ for every $-T<t<T$, or coincides with $\Sigma_i\times\{t\}$ for some $t$. The latter case implies $|S|_{g'_N}=|\Sigma_i|_{h_i}\geq A_0$, thus implies Claim 1. For the former case, we apply the above consequence of monotonicity formula for each slice $\Sigma_i\times[4j,4j+4]$, $|j|\leq\lfloor T/4\rfloor$, to obtain
		%\[|S|_{g'_N}\geq|S\cap\big(\Sigma_i\times(-T,T)\big)|_{g_i}\geq c_0\cdot2\lfloor\frac T4\rfloor>2A_0,\]
		%which contradicts the minimality of $S$ since one of the $\Sigma_i$ has area $A_0$.
		
		By the connectedness of $S$, we can assume without loss of generality that $S\subset M^-$. Since $M^-$ is topologically a sphere with disks removed, and $S$ does not bound a 3-ball in $M^-$, it follows that $S$ represents a nonzero element in $H_2(M^-,\ZZ)$. Now let $S'$ be a $g$-area minimizer in the homology class of $S$ in $H_2(M^-,\ZZ)$. We note that $S'$ must be a union of components of $\p M^-$. Indeed, any interior component of $S'$ is a smooth minimal surface, hence must be homologically trivial by the theorem's assumption. Thus we can decrease the area by removing this component, contradicting the minimizing property of $S'$. We then have $|S'|_g\geq A_0$. Combining what we obtained above, we have
		\[A_0 \leq |S'|_g
			\leq |S|_g
			= |S|_{g_N}
			\leq (1+\epsilon^2)|S|_{g'_N}
			= (1+\epsilon^2)\sys\pi_2(N,g'_N),\]
		which proves the claim. In the third inequality, we used the assumption $||g'_N-g_N||_{C^0(g_N)}\leq\epsilon^3$.
		% Detail: from the construction in Claim 2, we know that $S$ can be pushed into the interior of $M^-$ since the metric is warped product near $\p M^-$. So we can assume that $S$ is smooth in the interior of $M^-$.
	\end{proof}
	
	\begin{proof}[Proof of Claim 2] {\ }
		
		The proof involves smoothing $g_N$ near $\p M^\pm$ while preserving scalar curvature lower bounds. By symmetry, it suffices to perform the smoothing near $\p M^-$. We invoke the following gluing theorem in Brendle-Marques-Neves \cite{Brendle-Marques-Neves_2010}.
		
		\begin{theorem}[{\cite[Theorem 5]{Brendle-Marques-Neves_2010}}]\label{thm-sys:BMN}
			Let $M$ be a compact manifold with boundary $\p M$, and $g,\tilde g$ be two smooth metrics such that $g-\tilde g=0$ at each point on $\p M$. Moreover, assume that $H_{g}-H_{g'}>0$ at each point on $\p M$. Then given any number $\epsilon>0$ and any neighborhood $U$ of $\p M$, there is a smooth metric $\hat g$ on $M$ with the following property:
			
			(1) the scalar curvature of $\hat g$ satisfies $R_{\hat g}(x)\geq\min\{R_g(x),R_{\tilde g}(x)\}-\epsilon$ for every $x\in M$.
			
			(2) $\hat g$ agrees with $g$ outside $U$.
			
			(3) $\hat g$ agrees with $\tilde g$ in some neighborhood of $\p M$.
			
			(4) $||\hat g-g||_{C^0(g)}\leq\epsilon^3$.
		\end{theorem}
		
		Item (4) is not included in the original statement but follows directly from the proof. To achieve this, we choose the tensor $T$ to vanish outside a sufficiently small neighborhood of $\p M$ in \cite[p.189]{Brendle-Marques-Neves_2010}, then choose the parameter $\lambda$ to be sufficiently large in \cite[p.190]{Brendle-Marques-Neves_2010}.
		
		To apply the theorem, we need to perturb $g$ so that $\p M^-$ is strictly mean convex. Let $\delta$ be a small constant to be chosen. Choose any function $u\in C^\infty(\bar M)$ such that $||u||_{C^2(g)}\leq\delta$, $u|_{\p M}=0$ and $\frac{\p u}{\p\nu}>0$ on $\p M$. Set $g'=e^{2u}g$. The last condition of $u$ ensures that $\p M^-$ is mean convex with respect to $g'$. We have $||g'-g||_{C^0(g)}\leq\epsilon^3$ and $R_{g'}\geq\lambda-\epsilon/2$ when $\delta$ is sufficiently small.
		
		Next, we need to construct the metric $\tilde g$ which extends $g_i$ smoothly into the interior of $M^-$. We slightly enlarge $P_i$ to the cylinders $Q_i=\Sigma_i\times(-T-\delta,T]$, on which the metrics are still of warped product form $g_i=h_i+\varphi_i^2dt^2$. Let $\Phi_i:\Sigma_i\times(-T-\delta,-T]\to M^-$ be a regular smooth embedding, such that $\Phi_i$ maps $\Sigma_i\times\{-T\}$ identically to $\Sigma_i\subset\p M^-$, and $\Phi_i^*g'=g_i$ on $\Sigma_i\times\{-T\}$. Such a map can be constructed using normal exponential maps. Let $V_i\subset M^-$ be the image of $\Phi_i$. Thus, the metrics $\tilde g_i=(\Phi_i^{-1})^*g_i$ are defined in $V_i$ and coincides with $g'$ on $\Sigma_i\subset\p M^-$. Moreover, $\Sigma_i$ is totally geodesic with respect to $\tilde g_i$, hence $H_{g'}-H_{\tilde g_i}>0$ on $\Sigma_i$. Let $\tilde g$ be an arbitrary metric on $M^-$ that is equal to $\tilde g_i$ in a smaller neighborhood $U_i\subset V_i$ of $\Sigma_i$. Note that $R_{\tilde g}=R_{\tilde g_i}=(\Phi_i^{-1})^*R_{g_i}\geq\lambda$ in $U_i$. Now we apply Theorem \ref{thm-sys:BMN} in the neighborhood $\bigcup_i U_i\supset\p M^-$ and with $g'$ in place of $g$. We obtain a new metric $\hat g$ on $M^-$, such that $||\hat g-g'||_{C^0(g')}\leq\epsilon^3$ and $R_{\hat g}\geq\lambda-\epsilon$.
		
		Finally, set $g'_N$ to be equal to $\hat g$ on $M^-$, and equal to $g_i$ on $P_i$. Thus $g'_N$ satisfies the requirements of Claim 2 (on the side of $M^-$) once we specify a smooth structure across $\p M^-$ for which $g'_N$ is smooth. We express the space $X=M^-\cup_{\p M^-}(\bigcup P_i)$ as
		\[X=\Big(\text{int}(M^-)\sqcup\big(\bigcup Q_i\big)\Big)\Big/\!\sim,\]
		where $\sim$ is the equivalence relation $x\sim\Phi_i(x)$, $\forall x\in \Sigma_i\times(-T-\delta,-T)\subset Q_i$. Since the relation $\sim$ is given by diffeomorphism, a smooth structure on $X$ is naturally induced. It follows from Theorem \ref{thm-sys:BMN} (3) that $g'_N$ is smooth under this smooth structure. This completes the smoothing near $\p M^-$, and the smoothing near $\p M^+$ follows by symmetry. This proves the claim.
	\end{proof}
	
\end{proof}

\begin{remark}[{Deduction of Theorem \ref{thm-intro:sys_pi2} from Theorem \ref{thm-intro:local}}]\label{rmk-sys:decomp} {\ }
	
	Assume that Theorem \ref{thm-intro:local} holds. Let $M$ satisfy the assumptions of Theorem \ref{thm-intro:sys_pi2}. Following the idea in the introduction, we shall show that \eqref{eq-intro:sys_pi2} holds for any metric $g$ on $M$. We may assume without loss of generality that $g$ has positive scalar curvature; thus the topological classification \eqref{eq-sys:topology} is available. The argument is divided into several steps. We use $\#^k(S^2\times S^1)$ to denote the connected sum of $k$ copies of $S^2\times S^1$.
	
	\textit{Step 1:} we show that $M$ is finitely covered by $\#^k(S^2\times S^1)$ for some $k$. Suppose first that $M$ does not contain $S^2\times S^1$ prime factors, thus $M\cong(S^3/\Gamma_1)\#\cdots\#(S^3/\Gamma_m)$ for some finite groups $\Gamma_i$. By Kurosh's subgroup theorem \cite[Theorem VII.5.1]{Massey}, the kernel of the quotient map $\pi_1(M)\mapsto\Gamma_1\times\cdots\times\Gamma_m$ is a free group with finite index. Hence $M$ is finitely covered by a closed 3-manifold $\bar M$ with free fundamental group. By the classification \eqref{eq-sys:topology} again, we see that $\bar M\cong\#^k(S^2\times S^1)$ for some $k$. For the general case, we may write $M=\#^k(S^2\times S^1)\# N$, where $N$ is a connected sum of spherical space forms. From the special case, $N$ is covered by $\bar N=\#^m(S^2\times S^1)$ for some $m$. This implies that $M$ is convered by $\#^{k\deg(\bar N\to N)}(S^2\times S^1)\#\bar N$.
	
	\textit{Step 2:} we decompose $M$ into building pieces. Since the $\pi_2$\,--\,systole remains unchanged when passing to covering spaces, we may assume $M=\#^k(S^2\times S^1)$. By slightly perturbing the metric $g$, we assume that it is bumpy \cite{White_2017}. Let $\{\Sigma_1,\cdots,\Sigma_m\}$ be a maximal pairwise disjoint collection of stable minimal surfaces in $M$, such that each $\Sigma_i$ is nontrivial in $\pi_2(M)$. Such collection exists and is finite by the bumpiness hypothesis. We denote by $\{U_1,\cdots,U_n\}$ the set of connected components of $M\setminus\bigcup_{i\leq m}\Sigma_i$. The boundary of each $U_i$ consists of stable minimal spheres. Let $V_i$ be the closed manifold obtained by filling the boundaries of $U_i$ with 3-balls. Notice that $M$ can be recovered from $\bigsqcup_{i\leq n}V_i$ by performing 0-surgeries. Here, a 0-surgery means removing two disks and gluing the common sphere boundaries. If the two disks are in different (resp. the same) connected components of the manifold, then the 0-surgery is equivalent to a connected sum of the two components (resp. a connected sum of that component with $S^2\times S^1$). Therefore, from the 0-surgery we obtain
	\[M\cong \#^{m+1-n}(S^2\times S^1)\# V_1\#\cdots\# V_n.\]
	By the uniqueness of prime decomposition, each $U_i$ is diffeomorphic to either a 3-sphere with punctures or a $\#^s(S^2\times S^1)$ with punctures (for some $s\leq k$).
	
	We argue that the interior of each $U_i$ does not contain stable minimal surfaces that are nontrivial in $H_2(U_i,\ZZ)$. Suppose otherwise that there is a counterexample $\Sigma'\subset U_1$. In particular, $\Sigma'$ does not intersect any $\Sigma_i$. By the maximality of $\{\Sigma_i\}$, $\Sigma'$ must be trivial in $\pi_2(M)$. Then $\Sigma'$ must bound a simply-connected region (otherwise, any lift of $\Sigma'$ to the universal cover $\tM$ will be nontrivial in $H_2(\tM,\ZZ)$, contradicting the homotopy triviality of $\Sigma'$). By the uniqueness of prime decomposition, $\Sigma'$ bounds a 3-disk $D\subset M$. By the defining property of $\Sigma'$, the disk $D$ must contain one of the surfaces $\Sigma_i$. However, this implies that $\Sigma_i$ bounds a 3-disk, contradicting its homotopy triviality. This proves our claim.
	
	Note that, our claim implies that each $U_i$ must be a punctured $S^3$. Otherwise, one finds an area minimizer in $H_2(U_i,\p U_i,\ZZ)$ to obtain a contradiction.
	
	\textit{Step 3}: there exists a piece $U_1$ with at least three boundary components. Otherwise, all the $U_i$ are diffeomorphic to $S^2\times[0,1]$, which implies $M\cong S^2\times S^1$ and contradicts the assumption of Theorem \ref{thm-intro:sys_pi2}. Suppose $\Sigma_1,\cdots,\Sigma_l$ ($l\geq3$) are the boundaries of $U_1$. By the previous step, Theorem \ref{thm-intro:local} is applicable to $U_1$ and yields $\min_{i\leq l}|\Sigma_i|_g\cdot\min_{U_1}R_g\leq c_0$, where $c_0$ is the constant on the right hand side of \eqref{eq-intro:sys_pi2}. In particular, we have
	\[\sys\pi_2(M,g)\cdot\min_M R_g\leq\min_{i\leq l}|\Sigma_i|_g\cdot\min_{U_1}R_g\leq c_0.\]
\end{remark}

% New section. ------------------------------------------------------------

\vspace{12pt}

\noindent\textit{{Kai Xu,}}

\vspace{2pt}

\noindent\textit{{Department of Mathematics, Duke University, Durham, NC 27708, USA,}}

\vspace{2pt}

\noindent\textit{Email address: }\href{mailto:kx35@math.duke.edu}{kx35@math.duke.edu}.

\end{document}